\documentclass{amsart}
\usepackage[latin1]{inputenc}
\usepackage{graphicx}
\usepackage{subfigure}
\newtheorem{Teorema}{Theorem}
\newtheorem{Corolario}{Corollary}

\begin{document}

\bibliographystyle{alpha}
\title{Edge connectivity in graphs: an expansion theorem}
\author{Alvarez-Hamelin, José Ignacio \and Busch, Jorge Rodolfo}
\begin{abstract}
We show that if a graph is $k$-edge-connected,
and we adjoin to it another graph
satisfying a ``contracted diameter $\leq 2$'' condition,
with minimal degree $\geq k$,
and some natural hypothesis on the edges connecting one graph to the other,
the resulting graph is also $k$-edge-connected.
\end{abstract}
\maketitle

\section{Introduction}

Let $G$ be a simple graph ({\em i.e.}
a graph with no loops, no multiple edges) with vertex set $V(G)$ and edge set $E(G)$
(we follow in notation the book \cite{west:itgt}).
Given $A,B\subset V(G)$, $[A,B]$ is the set of edges
of the form $ab$,
joining a vertex $a\in A$ to a vertex $b\in B$.
As we consider edges without orientation,
$[A,B]=[B,A]$. Abusing of notation,
for $v\in V(G), A\subset V(G)$,
we write $[v,A]$ instead of $[\{v\},A]$.
The {\em degree} of a vertex $v\in V(G)$
is $\deg_G(v)\doteq |[v,V(G)]|$.
The {\em neighbourhood} of a vertex $v$, $N(v)$,
is the set of vertexes $w$ such that $vw \in E(G)$.
Given $A\subset V(G)$, $G(A)$
is the graph $G'$ such that $V(G')=A$ and $E(G')$
is the set of edges in $E(G)$ having both endpoints in $A$.
Given $v,w\in V(G)$, $d_G(v,w)$
is the distance in $G$ from $v$ to $w$,
that is the minimum length of a path from $v$ to $w$.
If $v\in V(G),A\subset V(G)$ we set
$d_G(v,A)\doteq \min_{w\in A} d_G(v,w)$.

An edge cut in $G$ is a set of edges $[S,\bar{S}]$,
where $S\subset V(G)$ is non void,
and $\bar{S}\doteq V(G)\setminus S$ 
is also assumed to be non void.

The edge-connectivity of $G$, $k'(G)$,
is the minimum cardinal of the cuts in $G$.
We say that $G$ is $k$-edge-connected if $k'(G)\geq k$.
Menger's theorem has as a consequence that given two vertices 
$v,w$ in $V(G)$,
if $G$ is $k$-edge-connected there are at least $k$-edge-disjoint paths
joining $v$ to $w$ (see \cite{west:itgt}, pp.153-169).

In this paper we address the following expansion problem: 
given a $k$-edge-connected graph $G_2$,
give conditions under which the result of adjoining
to $G_2$ a graph $G_1$ will be also $k$ edge-connected
(see Corollary \ref{conn} below).

\section{An expansion theorem }

Let $G$ be a simple graph.
Let $V_1\subset V(G),V_2\doteq V(G)\setminus V_1$, and set $G_1=G(V_1),G_2=G(V_2)$. We assume in the sequel that $V_1$ and $V_2$ are non void.
We define, for $x,y\in V_1$, the {\em contracted distance}
\[
\delta(x,y)\doteq \min\{d_{G_1}(x,y),d_G(x,V_2)+d_G(y,V_2)\}\\
\]
and for $x\in V,y\in V_2$
\[
\delta(x,y)=\delta(y,x)\doteq d_G(x,V_2)
\]
If $x\in V$ and $A\subset V$,
we set $\delta(x,A)\doteq \min_{a\in A} \delta(x,a)$.

Notice that with these definitions,
if $\delta(x,y)=2$ for some $x,y\in V$, then there exists $z\in V$ 
such that $\delta(x,z)=\delta(z,y)=1$.

We shall also use the notations 
\begin{eqnarray*}
\partial^j V_1&\doteq&\{x\in V_1: |[x,V_2]|\geq j\}\\
i^j V_1&\doteq& \{x\in V_1: |[x,V_2]|< j\}
\end{eqnarray*}

\begin{figure}[hb]
  \subfigure[]{
    \includegraphics[width=0.4\textwidth]{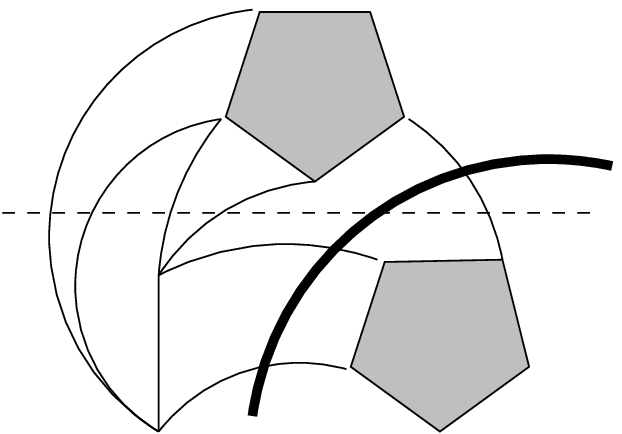}
  }
  \subfigure[]{
    \includegraphics[width=0.4\textwidth]{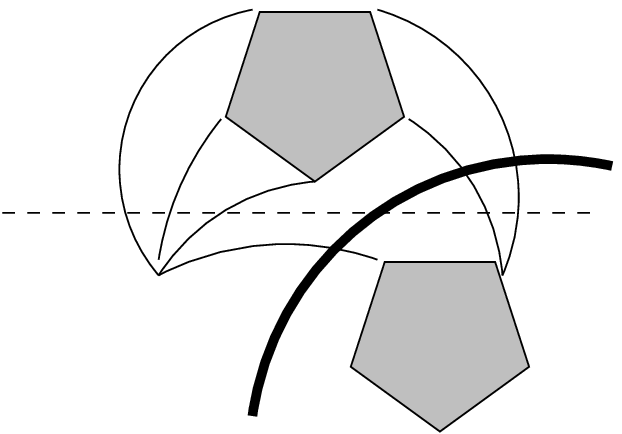}
  }
  \subfigure[]{
    \includegraphics[width=0.4\textwidth]{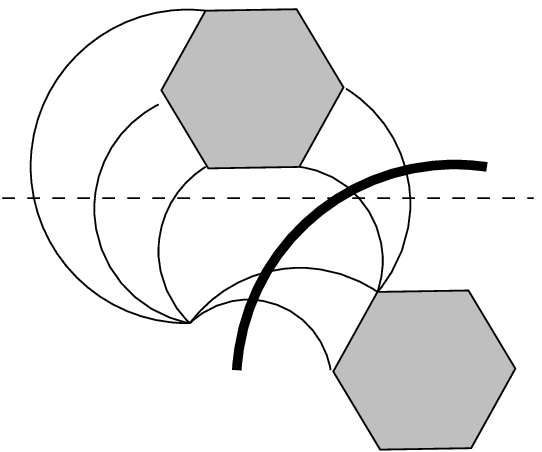}
  }
\hspace{0.1\textwidth}
  \subfigure[]{
    \includegraphics[width=0.35\textwidth]{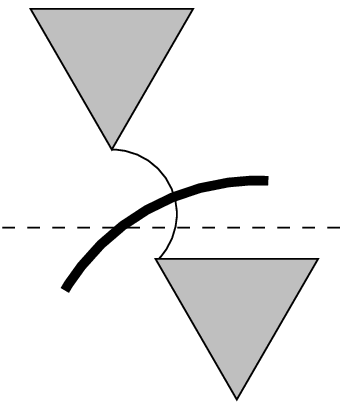}
  }
  \caption{\label{thm.ex} 
    {\bf Conventions:}
    1.Filled polygons represent cliques, and curved arcs represent edges.
    2.The dotted line separates 
    $G_2$ (the upper graph) from 
    $G_1$ (the lower graph).
    3. The widest arc shows the cut $[S,\bar{S}]$.
    {\bf Descriptions:} 
(a) Here $|[S,\bar{S}]|=3<k=4$, $|S\cap V_1|=2$,
$|\bar{S}|=5$, $\Phi=3$, $S\cap V_1=\partial^2 V_1$.
(b) Here $|[S,\bar{S}]|=3<k=4$, $|S\cap V_1|=1$,
$|\bar{S}|=5$, $\Phi=3$, $S\cap V_1\not =\partial^2 V_1$. 
(c) Here $|[S,\bar{S}]|=4<k=5$, $|S\cap V_1|=1$,
$|\bar{S}|=6$, $\Phi=3$, $S\cap V_1\not =\partial^2 V_1$. 
(d) Here $|[S,\bar{S}]|=1<k=2$, $|S\cap V_1|=0$,
$|\bar{S}|=3$, $\Phi=1$, $S\cap V_1 =\partial^2 V_1=\emptyset$. 
}
\end{figure}

Under these settings, we consider also
\[
\Phi\doteq\sum_{x\in V_1}\min\{\max\{1,|[x,i^2 V_1]|\},|[x,V_2]|\}
\]
In this general framework, we have
\begin{Teorema}
\label{cuts}
If $\max_{x,y\in V} \delta(x,y)\leq 2$
({\em i.e.} the {\em contracted diameter} of $V$ is $\leq 2$), 
$[S,\bar{S}]$ is an edge cut in $G$ such that
$V_2\subset S$,
and $k\doteq \min_{x\in V_1} \deg_G(v) >|[S,\bar{S}]|$,
then
\begin{enumerate}
\item
$\exists \bar{s}\in \bar{S}:\delta(\bar{s},S)=2$.
\item
$\forall s\in S: \delta(s,\bar{S})=1$.
\item
$|S\cap V_1|<|[S,\bar{S}]|<k<|\bar{S}|$.
\item
$S\cap V_1\subset \partial^2 V_1, \bar{S}\supset i^2 V_1$.
\item
$
\Phi \leq |[S,\bar{S}]|
$.
\end{enumerate}
\end{Teorema}
(See the examples in Figure \ref{thm.ex}.)

Proof.
\begin{enumerate}
\item
Arguing by contradiction, suppose that 
for any $\bar{s}\in \bar{S}$: $\delta(\bar{s},S)=1$.
Let $\bar{s}\in\bar{S}$.
Then we have $k_1$ edges $\bar{s}s_i,1\leq i\leq k_1$
with $s_i\in S$ and (eventually) $k_2$ edges $\bar{s}\bar{s}_j$,  
$\bar{s}_j\in \bar{S}$.
But each $\bar{s}_j$ satisfies $\delta(\bar{s}_j,S)=1$,
thus we have $k_2$ new edges (here we used that $G$ is simple,
because we assumed that the vertices $\bar{s}_j$ are different)
$\bar{s}_js'_j$, with $s'_j\in S$,
whence \[|[S,\bar{S}]|\geq k_1+k_2=\deg_G(\bar{s})\geq k\]
which contradicts our hypothesis.
\item
Let $\bar{s}_0\in \bar{S}$ be such that $\delta(\bar{s}_0,S)=2$.
Then for each $s\in S$, as  $\delta(\bar{s}_0,s)=2$,
there exists $\bar{s}'$ such that 
$\delta(\bar{s}_0,\bar{s}')=\delta(\bar{s}',s)=1$.
But, again, as $\delta(\bar{s}_0,S)=2$, it follows that $\bar{s}'\in \bar{S}$,
hence $\delta(s,\bar{S})=1$.
\item
By the previous point, 
we have for each $s\in S\cap V_1$ 
some edge in $[S,\bar{S}]$ incident in $s$,
and for some $v\in V_2$ we have also
some edge in $[S,\bar{S}]$ incident in $v$,
thus
\[
|S\cap V_1|+1\leq |[S,\bar{S}]|
\]
On the other hand, 
if $\bar{s}\in \bar{S}$ satisfies $\delta(\bar{s},S)=2$
(such $\bar{s}$ exists by our first point),
then $N(\bar{s})\subset \bar{S}$ 
(recall that $N(\bar{s})$ is the neighbourhood of $\bar{s}$),
whence
\[
|\bar{S}|\geq 1+|N(\bar{s})|\geq 1+k
\]
and our statement follows.
\item
Let $s\in S\cap V_1$.
By our second point, and using
again that there is at least one edge in $[\bar{S},V_2]$,
we have
\begin{eqnarray*}
|N(s)\cap V_1|+|[s,\bar{S}]|&\leq& |[S,\bar{S}]|-1\\
&<&\deg_G(s)-1
\end{eqnarray*}
and the first of our statements follows if we notice that
\[
\deg_G(s)=
|N(s)\cap V_1|+|[s,\bar{S}]|+|[s,V_2]|
\]
Now, $\bar{S}=V_1\setminus S\cap V_1$,
and our second statement follows immediately.
\item
By our previous points,
if $s\in S\cap V_1$ then
\[
|[s,\bar{S}]|\geq \max\{1,|[s,i^2V_1]|\}
\]
and of course for $\bar{s}\in \bar{S}$, $|[\bar{s},S]|\geq |[\bar{s},V_2]|$,
thus
\begin{eqnarray*}
|[S,\bar{S}]|&=&|[S\cap V_1,\bar{S}]|+|[\bar{S},V_2]|\\
&\geq&\sum_{s\in S\cap V_1}\max\{1,|[s,i^2V_1]|\}+\sum_{\bar{s}\in \bar{S}} |[\bar{s},V_2]|\\
&\geq&
\Phi
\end{eqnarray*}
\end{enumerate}



\begin{figure}[hb]
  \subfigure[]{
    \includegraphics[width=0.4\textwidth]{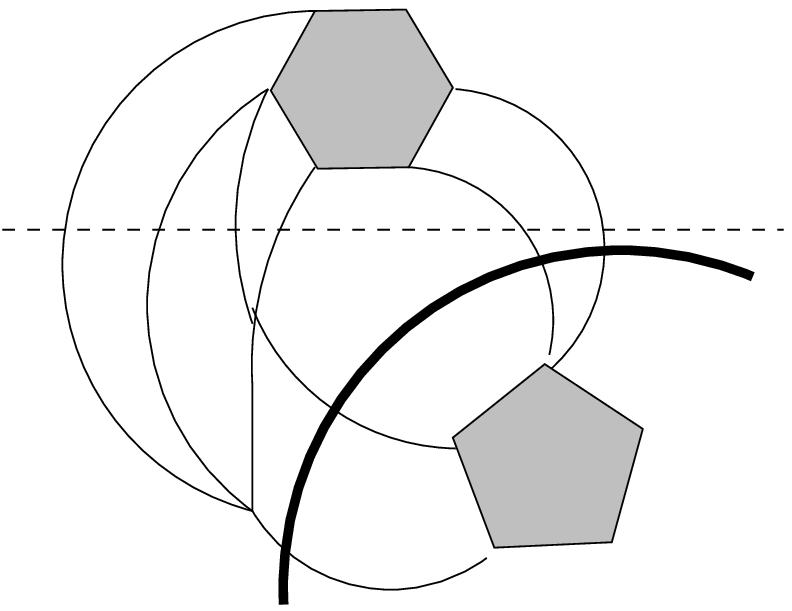}
  }
  \subfigure[]{
    \includegraphics[width=0.4\textwidth]{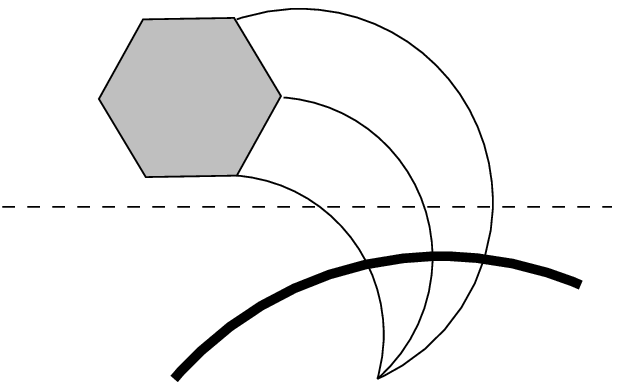}
  }
  \caption{\label{crll.ex} 
    {\bf Conventions:}
    1.Filled polygons represent cliques, and curved arcs represent edges.
    2.The dotted line separates 
    $G_2$ (the upper graph) from 
    $G_1$ (the lower graph).
    3. The widest arc shows a minimal cut $[S,\bar{S}]$.
    {\bf Descriptions:} 
(a) Here $|[S,\bar{S}]|=k=4$, $\Phi=4$, $|\partial^1 V_1|=3$.
(b) Here $|[S,\bar{S}]|=k=3$, $\Phi=1$, $|\partial^1 V_1|=1$,
$V_1=\partial^1 V_1$. 
This example shows that Corollary \ref{conn} includes an 
edge-connectivity version
of the Expansion Lemma in \cite{west:itgt}, Lemma 4.2.3.
}
\end{figure}

\begin{Corolario}
\label{conn}
Assume that 
\begin{enumerate}
\item
$\deg_G(x)\geq k,x\in V(G)$
\item
$G_2$ is $k$-edge connected
\item
$\max_{x,y\in V} \delta(x,y)\leq 2$
\end{enumerate}
Then any of the following 
\begin{enumerate}
\item
$\Phi\geq k$
\item
$|\partial^1 V_1|\geq k$ 
\item
$V_1=\partial^1 V_1$
\end{enumerate}
implies that  $G$ is $k$-edge-connected.
\end{Corolario}
(See the examples in Figure \ref{crll.ex}.)

Proof.
Let $[S,\bar{S}]$ be any cut in $G$. 
We shall show that,
under the listed hypotheses and any of
the alternatives, $|[S,\bar{S}]|\geq k$.

If $S\cap V_2\not =\emptyset$ and $\bar{S}\cap V_2\not = \emptyset$,
then, as
\[
[S\cap V_2,\bar{S}\cap V_2]\subset[S,\bar{S}]
\]
is a cut in $G_2$, which we assumed to be $k$-edge connected,
we obtain $|[S,\bar{S}]|\geq k$.

Without loss of generality, we assume in the sequel that
$V_2\subset S$. We argue by contradiction
assuming that there exists some $S$ such that $|[S,\bar{S}]|<k$,
so that we are under the hypothesis of Theorem \ref{cuts}.

The first of our alternative hypothesis contradicts
the last of the conclusions of Theorem \ref{cuts}.

When $x\in \partial^1V_1$, 
\[
\min\{\max\{1,|[x,i^2 V_1]|\},|[x,V_2]|\}\geq 1\] 
so that we have
$
|\partial^1 V_1| \leq\Phi
$
{\em i.e.} the second of our alternative hypothesis implies the first one.

To finish our proof,
notice that if $V_1=\partial^1 V_1$, as $\bar{S}\subset V_1$,
we have $\delta(\bar{s},S)=1$ for any $\bar{s}\in \bar{S}$,
contradicting the first of the conclusions in Theorem \ref{cuts}.

\section{Final remarks}

Corollary \ref{conn} is 
related to a well known theorem of Plesník 
(see \cite{plesnik:cgoagd}, Theorem 6),
which states that in a simple graph
of diameter $2$ the edge connectivity
is equal to the minimum degree.
\bibliography{arc}

\end{document}